%%%%%%%%%%%%%%%%  Macros specific to this document %%%%%%%%%%%%%%%%%%%

\def\label#1{\global\edef#1{\the\sectnumber.\the\thenumber}}

\newcount\myequno
\myequno=0
\def\myeqno{\global\advance\myequno by 1 \eqno(\the\myequno)}
\def\eqlabel#1{\global\edef#1{\the\myequno}}

\def\MMProp#1{\Prop{\sl #1}}
\def\MMLemma#1{\Lemma{\sl #1}}

\def\cA{{\cal A}}

\def\cH{{\cal H}}

\def\be{{\bf e}}
\def\bn{{\bf n}}

\def\bx{{\bf x}}
\def\by{{\bf y}}
\def\bz{{\bf z}}

\def\bU{{\bf U}}

\def\tf{{\tilde f}}

\def\tx{{\tilde x}}

\def\tz{{\tilde z}}
\def\tE{{\tilde E}}

\def\tZ{{\tilde Z}}
\def\tepsilon{{\tilde \epsilon}}

\def\txi{{\tilde \xi}}

\def\Cstar{{\bf C}^*}

\def\cA{{\cal A}}

\def\Delj{\Delta(n_j)}

\def\vol{{\rm vol}}
\def\per{{\rm per}}

\def\Expect#1#2{{\bf E}_{#1}\big(#2\big)}

\def\IndexGoal{E(\bn, \delta)}

\def\BKK{B\!K\!K}

\def\Gbig{\Big}
\def\Gbigl{\Bigl}
\def\Gbigr{\Bigr}

\font\elevenbf=cmbx10 scaled \magstephalf  % boldface for titles 
\font\twelvebf=cmbx10 scaled 1315          % large boldface for large titles 

\def\Re{I\!\!R}
\def\bull{\vrule height 1.5ex width 1.0ex depth -.1ex}

\newcount\sectnumber
\newcount\subsectnumber
\newcount\eqnnumber
\newcount\notenumber
\newcount\thenumber
\newcount\teoremnumber
\def\eqnno#1{\global\advance\eqnnumber by 1 $$#1\eqno(\the\eqnnumber)$$}

\def\note#1{\global\advance\notenumber by 1 
          \footnote{$^{(\the\notenumber)}$}{#1}}
\def\Ass{\noindent \global\advance\thenumber by 1 
           {\bf Assumption $\bf \the\sectnumber$.$\bf \the\thenumber$:}\enspace}

\def\Not{\noindent \global\advance\thenumber by 1
           {\bf Notation $\bf \the\sectnumber$.$\bf \the\thenumber$:}\enspace}

\def\Def{\noindent \global\advance\thenumber by 1
           {\bf Definition $\bf \the\sectnumber$.$\bf \the\thenumber$:}\enspace}

\def\Rem{\noindent \global\advance\thenumber by 1
           {\bf Remark $\bf \the\sectnumber$.$\bf \the\thenumber$:}\enspace}

\def\Th{\noindent \global\advance\teoremnumber by 1
           {\bf Theorem $\bf \the\teoremnumber$:}\enspace}
\teoremnumber = 0

\def\Cor{\noindent \global\advance\thenumber by 1
           {\bf Corollary $\bf \the\sectnumber$.$\bf \the\thenumber$:}\enspace}

\def\Ex{\noindent \global\advance\thenumber by 1
          {\bf Example $\bf \the\sectnumber$.$\bf \the\thenumber$:}\enspace}

\def\Prop{\noindent \global\advance\thenumber by 1
            {\bf Proposition $\bf \the\sectnumber$.$\bf \the\thenumber$:}\enspace}
\def\AppBProp{\noindent \global\advance\thenumber by 1
            {\bf Proposition B.$\bf \the\thenumber$:}\enspace}

\def\Claim{\noindent \global\advance\thenumber by 1
            {\bf Claim $\bf \the\sectnumber$.$\bf \the\thenumber$:}\enspace}

\def\Lemma{\noindent\global\advance\thenumber by 1
             {\bf Lemma $\bf \the\sectnumber$.$\bf \the\thenumber$:}\enspace}
\def\AppALemma{\noindent\global\advance\thenumber by 1
             {\bf Lemma A.$\bf \the\thenumber$:}\enspace}
\def\endlemma{\bigskip\rm}
  
\def\endproof{\enspace\bull\bigskip}

% The following were added by AM, in response to a desire for
% macros that did not constrain theorem and lemma numbers unduly

%%%%%%%%%

\def\Section#1{\bigskip
\eqnnumber=0
\thenumber=0
\subsectnumber=0
               \global\advance\sectnumber by 1
               {\noindent{\elevenbf$\bf \the\sectnumber$.\enspace#1}}
                 \nobreak}

\def\Subsection#1{\medskip
		\global\advance\subsectnumber by 1
            {\noindent $\quad$\it$\bf \the\sectnumber$.$\bf \the\subsectnumber$.\enspace#1}
                 \smallskip}

%%%%%%%%%%%%%%%%%% Start of MSRI Form %%%%%%%%%%%%%%%%%%%%%%%%%

\magnification=1200 
\vsize=8truein 
\hsize=6truein 
\hoffset=18truept
\voffset=24truept

\font\smc=cmcsc10
\font\smallsmc=cmcsc8
\font\smallrm=cmcsc8

\headline={\ifodd\pageno \ifnum\pageno>1 \smallrm \hfil 
EXPECTED NUMBER OF REAL ROOTS                % running head for right-hand page is title in caps
\hfil\folio \else\hfill\fi \else \smallrm \folio \hfill
ANDREW MCLENNAN   % running head for left-hand page is authors in caps
\hfill\fi} \footline={\hss}   % footline is blank

\vglue .5in

\centerline{\bf
THE EXPECTED NUMBER OF REAL ROOTS}
\centerline{ \bf OF A MULTIHOMOGENEOUS SYSTEM 
}
\centerline{\bf OF
POLYNOMIAL EQUATIONS                 %use capital letters for title
}
\vglue .5in
\centerline{\smc        %cap and lowercase for author
Andrew McLennan
}
\vglue .5in
\narrower{\noindent
{\smc Abstract.}
The methods of Shub and Smale [SS93] are extended to the
class of multihomogeneous systems of polynomial equations, yielding
Theorem 1, which is a formula expressing the mean (with respect to a
particular distribution on the space of coefficient vectors) number of
real roots as a multiple of the mean absolute value of the determinant
of a random matrix. Theorem 2 derives closed form expressions for the
mean in special cases that include: (a) Shub and Smale's result that
the expected number of real roots of the general homogeneous system is
the square root of the generic number of complex roots given by
Bezout's theorem; (b) Rojas' [Roj96] characterization of the mean
number of real roots of an ``unmixed'' multihomogeneous system.
Theorem 3 gives upper and lower bounds for the mean number of roots,
where the lower bound is the square root of the generic number of
complex roots, as determined by Bernstein's [Ber75] theorem.  These
bounds are derived by induction from recursive inequalities given in
Theorem 4.
}
\vglue .5in
\footnote{}{\hskip-\parindent
Research at MSRI is supported in part by NSF grant DMS-9701755.
I have benefited from numerous discussions with Maurice Rojas, and I
am grateful for comments by seminar participants at the University of
Minnesota, the AMS Meeting at Temple University in April 1998, and the
conference on polynomial system solving at the Mathematical Sciences
Research Institute in September 1998.  }

\vfill \eject

\centerline{\twelvebf The Expected Number of Real Roots of a}
\centerline{\twelvebf Multihomogeneous System of Polynomial Equations}

\bigskip

\vskip .1cm

\Section{Introduction}

\medskip
The study of the distribution of real roots of a polynomial with random coefficients, which traces back 
at least to [BP32], has recently been developed in the direction of multivariate systems.
(This literature is
ably surveyed, and extended, by Edelman and Kostlan [EK95].)
Kostlan [Kos93] shows that, for a homogeneous polynomial equation of 
degree $d$ in $n+1$ variables, a particular inner product on the space of coefficient vectors is
distinguished by invariance under the natural action of $O(n+1)$ and orthogonality of monomials. He goes on to show that, for the system of $n$ such equations, when the coefficient vectors for the various equations are independent random variables, with each one distributed according to the central normal distribution associated with this inner product, the mean number of projective roots in $n$-dimensional real projective space is $d^{n/2}$, which is the square root of the generic number of complex roots given by Bezout's theorem.
Shub and Smale [SS93] extend this result to the general homogeneous system of $n$ homogeneous polynomial
equations of degrees $d_1, \ldots, d_n$, showing that the mean is
$\sqrt{\prod_i d_i}$, which is again the square root of the Bezout number.
Rojas [Roj96] studies unmixed\note{Sparse systems of polynomial equations are described in Section 2.  Roughly, such a system is unmixed if all polynomials have the same collection of monomials with nonzero coefficients, and otherwise it is mixed.} systems of multihomogeneous equations,
arriving at a closed form formula for the mean number of roots in the cartesian product of projective spaces that is the natural root space for such systems.  

This paper studies the more general case of mixed multihomogeneous systems.  Theorem 1 
is a formula expressing the mean number of real roots of a random multihomogeneous system  as the product of the mean absolute value of the
determinant of a random matrix times an expression composed of
evaluations of Euler's function $\Gamma$ at multiples of $1/2$. 
Theorem 2, which is a corollary, gives a closed form formula for this mean, for a smaller
class of systems that includes both the general homogeneous system and
the unmixed systems as special cases, so that the results of [SS93] and [Roj96] described above are corollaries. 
Theorem 3 generalizes the ``square root of the Bezout number'' result
by giving upper and lower bounds on the mean number of roots, where
the lower bound is the square root of the maximal number of roots for
the associated ``demultihomogenized'' system, as given by Bernshtein's [Ber75] extension of Bezout's theorem to sparse systems of polynomial equations.  These bounds follow from
recursive inequalities given in Theorem 4. 

The author's interest in this topic is motivated in part by concepts
of noncooperative game theory\note{This is not the place to give a
general introduction to noncooperative game theory; Fudenberg and
Tirole (1991) is a standard text.  For the internal logic of this
paper the description of {\it quasiequilibrium\/} (Section 2) is
sufficient. For the connection between this notion and the standard
concepts of {\it Nash equilibrium\/} and {\it totally mixed Nash
equilibrium\/} see [MM97, McL97].}.  The concept of a totally mixed Nash equilibrium for a normal form game amounts to a root, all of whose components must be positive, of particular sort of multihomogeneous system.  
McLennan and McKelvey [MM97] give
a method for constructing normal form games that have as many regular (real)
totally mixed Nash equilibria as are permitted by Bernshtein's
theorem.  The conceptual import of this result is that the maximal
number of Nash equilibria is large, at least compared to most game
theorists' prior intuition.  Games that have the maximal number of
equilibria are thought to be very atypical, and there arises the
question of whether the set of equilibria is not only potentially
large, but also large on average.  McLennan [McL97]
investigates the application, to this problem, of the results
developed here, using Theorem 3 to show that the mean number of Nash
equilibria can grow exponentially with various measures of the size of
the game.  Among other things, this analysis involves the extension of
our work here to systems consisting of a multihomogeneous system of
the sort studied here to which additional multihomogeneous polynomial
inequalities have been appended, with the generalized formula being
the one given here times a factor that may be regarded as the
``probability'' that a root of the system of equations also satisfies
the inequalities. 

In connection with speculation concerning whether analogues of Theorem
3 might hold for more general classes of sparse systems than the
multihomogeneous ones, we recommend [Roj], which gives an extension to
general sparse systems of the model of a random system studied here,
and which presents results and conjectures along these lines.  It is
interesting to note that multihomogeneous systems are potentially
special insofar as they can have as many real regular roots as are
permitted by Bernshtein's theorem.  (This is proved in [McL98] by
pointing out that the argument in [MM97], which establishes this claim
for the systems arising in game theory, is actually valid for any
multihomogeneous system.)

The proof of Theorem 1 parallels the analysis in [SS93] and [BCS98] rather closely, and is thus a descendant of the methods of [Kac43].  
The {\it incidence variety\/} is the set of coefficient vector-root pairs.  It is a submanifold of the cartesian product of the space of coefficient vectors and the root space, and the projection of it onto the root space is a fibration.  The roots of the system at a particular coefficient vector are the preimages of the projection of the incidence variety onto the space of coefficient vectors, and an integral formula [SS93, p.~273;  BCS98, p.~240] is used to reexpress the mean number of roots as a double integral, where the outer integral is over the root space and the inner integral is over the fibre of the projection onto the roots space at the root in question.  
Invariance is used to show that the inner integral does not depend on this root, so that the double integral is the volume of the root space times the inner integral, evaluated at a point in the root space which may be chosen at whim.  For a particular choice it is possible to simplify the inner integral by transforming variables in a way that eliminates variables that do not enter the integrand, and from this Theorem 1 emerges.

The algorithms used by [MM97] to compute maximal numbers of Nash
equilibria are based on recursive formulas for the Bernshtein number
that extend directly to general multihomogeneous systems.  Below (see
also [McL98]) we describe how these formulas can be seen as the
consequence of expressing the Bernshtein number for such a system as
the permanent (e.g.~[Ego96]) of a matrix, after which the recursions
are obtained by expanding along a row or column.  In investigating whether the mean number of
real roots is greater than the square root of the Bernshtein number,
as asserted by Theorem 3, it is natural to guess that the squares of
the mean numbers of real roots obey the corresponding recursive
inequalities, which is the assertion of Theorem 4, since then Theorem
3 follows from induction.  Using Theorem 1, Proposition 7.1 restates
these inequalities as recursive inequalities for the mean absoute
values of the determinants of certain random matrices.  The proof of
Proposition 7.1 is, perhaps, rather surprising insofar as it depends
on properties of normal random variables that seem quite distant from
the geometric starting point of these investigations.

The remainder has the following organization.  Section 2 describes
multihomogeneous systems as a certain type of sparse system.  Section
3 specifies an inner product on the space of coefficient vectors of a
multihomogeneous equation that is uniquely characterized by invariance and orthogonality of monomials.  The central normal distribution with respect
to this inner product is our model of a random equation, and our random systems have the coefficient vectors of the various equations distributed independently according to these distributions.  Section 4
states Theorem 1, and in Section 5 we discuss those systems for which
it is possible to reduce the formula in Theorem 1 either to closed
form or to an expression involving the formula applied to smaller
systems.  Section 6 defines mixed volume, states Bernshtein's theorem
precisely, and shows how the generic number of complex roots of a
multihomogeneous system may be computed recursively.  Section 7 proves
Theorems 3 and 4, and presents a result giving upper and lower bounds
for the mean absolute value of the determinant of a random matrix.
Sections 8--11 present the proof of Theorem 1.

\vskip .1cm

\medskip
\Section{Multihomogeneous Systems}

\medskip
In stating Bernshtein's theorem we will need to consider general
sparse systems, so we describe multihomogeneous systems as a
specialization of this concept.  A {\it sparse system\/} of $n$
polynomial equations in $\ell$ variables is $$f(\bx) = (f_1(\bx),
\ldots, f_n(\bx)) = 0,
$$ where $\bx\,=\,(x_1,\ldots,x_{\ell})\,$ and, for each $i = 1,
\ldots, n$, there is a nonempty finite $\cA_i
\subset {\bf N}^{\ell}$ such that $f_i (\bx) = \sum_{a \in \cA_i}f_{ia}
\bx^a$ for some system of coefficients $f_{ia}$. (Here 
$\bx^a$ denotes the monomial $x_1^{a_1}x_2^{a_2}\cdots
x_{\ell}^{a_{\ell}}$.)  The general approach of the theory of sparse
systems is to hold the $n$-tuple of {\it supports\/} $(\cA_1, \ldots,
\cA_n)$ fixed while treating the coefficients $f_{ia}$ as variables,
for instance in the sense of studying properties that are generic in
the space of vectors of coefficients.  Such a system is said to be
{\it unmixed\/} if $\cA_1 = \ldots = \cA_n$; otherwise it is {\it
mixed}.  Identifying a polynomial with its vector of coefficients, we
regard $\cH_i := \Re^{\cA_i} \setminus \{0\}$ as the space of
polynomials with real coefficients whose supports are nonempty subsets
of $\cA_i$.  Let $$\cH := \cH_1 \times \ldots \times \cH_n.$$

The system is {\it multihomogeneous\/} if the variables in $\bx$ are
divided into $k$ groups, so that $\bx = (\by_1, \ldots, \by_k)$ where
$\by_j = (y_{j0},y_{j1}, \ldots, y_{jn_j})$, and each equation is
homogeneous of degree $\delta_{ij}$ as a function of $\by_j$, for any
given values of the other variables $(\by_1, \ldots, \by_{j-1},
\by_{j+1}, \ldots, \by_k)$.  More precisely, we require that
there are nonnegative integers $\delta_{ij}$ ($i = 1, \ldots, n$, $j =
1, \ldots, k$) such that $$\cA_i = \cA_{i1} \times \ldots \times
\cA_{ik}, \quad \hbox{where}
\quad \cA_{ij} = \{\, \alpha \in {\bf N}^{n_j+1} : \alpha_0 + \alpha_1 +
\ldots + \alpha_{n_j} = \delta_{ij} \,\}.$$
When $f_i$ is multihomogeneous, the truth value of the proposition
`$f_i(\bx) = 0$' is unaffected if each block of variables is
multiplied by a nonzero scalar, so that, in effect, there are $\ell -
k$ degrees of freedom.  We work only with systems that are, in this
sense, exactly determined: $\ell = n + k$, so that $n_1 + \ldots +
n_k = n.$ An instance of the type of system studied here is specified
by the vector $\bn$ and the $n \times k$ matrix $\delta :=
(\delta_{ij})$.

Four particular types of multihomogeneous system figure in our
discussion:
\smallskip
\item{(a)}  When $k = 1$ we have the general homogeneous system,
for which the problem studied here was analyzed in [SS93].  In
inductive constructions it will be convenient to allow the numbers of
variables in some blocks to be zero, and we will use the phrase
`general homogeneous system' to describe any multihomogeneous system
with $n_j = n$ for some $j$, in which case we must have $n_h = 0$ for
all $h \ne j$.
\item{(b)}  The
{\it unmixed\/} multihomogeneous systems studied in [Roj96] are
described by the condition that all equations have the same support:
there are integers $e_1, \ldots, e_k$ such that $$\delta_{1j} = \ldots
= \delta_{nj} = e_j \quad (j = 1,
\ldots, k).$$
\item{(c)} Generalizing (a) and (b) are the systems for which there
are numbers $d_1, \ldots, d_n$ and $e_1, \ldots, e_k$ such that
$\delta_{ij} = d_ie_j$ for all $i$ and $j$.
\item{(d)} The systems
arising, in game theory, from the concept of quasiequilibrium ([MM97])
of a finite normal form game, have, for each $j = 1, \ldots,k$, $n_j$
equations that are homogeneous of degree one in $\by_h$ for all $h \ne
j$, and are homogeneous of degree zero in $\by_j$.  Formally these
systems can be characterized as follows: $$\delta_{ij} = \cases{ 0 &
if $q(i) = j$, \cr
\noalign{\smallskip}
  1  &  otherwise,     \cr
}\myeqno \eqlabel{\NashEquilibrium}
$$
where $q : \{1, \ldots, n\} \to \{1,
\ldots, k\}$ is the function defined implicitly by the inequality
$$n_1 + \ldots + n_{q(i)-1} < i \le n_1 + \ldots + n_{q(i)}.$$

\vskip .1cm

\goodbreak
\medskip
\vfill \eject
\Section{An Invariant Inner Product}
\nobreak
\medskip
Fix a pair $(\bn,\delta)$.  Since $n_1 + \ldots + n_k = n$, we may
index the components of an exponent vector $a \in {\bf N}^{n+k}$ by
the pairs $(j,h)$ for $j = 1, \ldots, k$ and $h = 0, \ldots, n_j$.
For such an $a$ let $$\eqalign{
\eta  &  (a)   := { a_{10}! \cdot \ldots \cdot a_{1n_1}! \over
(a_{10} + \ldots + a_{1n_1})!} \cdot
\ldots \cdot { a_{k0}! \cdot \ldots \cdot a_{kn_k}! 
\over (a_{k0} + \ldots + a_{kn_k})!} \cr
&  \strut \cr
&
= \left( { a_{10} + \ldots + a_{1n_1} \atop a_{10}, \ldots, a_{1n_1} }
\right)^{-1} \cdot \ldots \cdot \left( { a_{k0} + \ldots + a_{kn_k} 
\atop a_{k0}, \ldots,
a_{kn_k} } \right)^{-1}. \cr
}$$ 
We endow each $\cH_i$ with the inner product $$\langle f_i, f'_i
\rangle_i :=
\sum_{a \in \cA_i} \eta(a) f_{ia}f'_{ia}.$$
Let $\|\cdot\|_i$ be the norm derived from $\langle \cdot, \cdot
\rangle_i$.  

Consider the product group $$G := O(n_1 + 1) \times \ldots \times
O(n_k + 1).$$ There is the obvious component-wise action of $G$ on
$\Re^{n_1+1} \times \ldots \times \Re^{n_k+1}$, and for $f_i \in
\cH_i$ and $O \in G$, $f_i \circ O^{-1}$ is easily seen to be a
polynomial function that is multihomogeneous for the same numbers
$\delta_{ij}$, so $f_i
\circ O^{-1}$ is an element of $\cH_i$.  Thus the formula $Of_i := f_i
\circ O^{-1}$ defines an action from the left of $G$ on $\cH_i$.  The following generalizes [Kos93, Th.~4.2], which is the case $k = 1$.  

\medskip
\MMLemma{ 
\label{\UniqueInvariance}
The inner product (4) is the unique (up to multiplication by a scalar)
inner product on $\cH_i$ that is invariant under the action of $G$ and
with respect to which the monomials are pairwise orthogonal.
}
\endlemma

\medskip
\noindent{\bf Proof:}  
Let $\langle f_i, f'_i
\rangle_i^* =
\sum_{a \in \cA_i} \eta^*(a) f_{ia}f'_{ia}$ be an invariant inner
product with all monomials orthogonal.  We wish to show that
$\eta^*(a)/\eta^*(a') = \eta(a)/\eta(a')$ for all $a, a' \in \cA_i$.
Fixing arbitrary $a \in \cA_i$ and $j = 1, \ldots, k$, it suffices to
establish that this formula holds for those $a' \in \cA_i$ with
$a_{\ell h} = a'_{\ell h}$ whenever $\ell \ne j$, and this follows
from [Kos93, Th.~4.2] applied to the subspace of $\cH_i$ spanned by
such $a'$.

To see that $\langle \cdot , \cdot \rangle_i$ is invariant under the
action of $G$ observe that, by [Kos93, Th.~4.2], it is invariant under
the action of any group element $g$ with only one component
$g_j$ different from the identity in $O(n_j + 1)$, and that such group
elements generate $G$. \endproof

Following [Kos93, EK95, Roj96], in our model of a random multihomogeneous
system  the coefficient vectors of the various
equations are statistically independent, with the coefficient vector of
the $i^{\rm th}$ equation centrally normally distributed in $\cH_i$ relative to
$\langle \cdot, \cdot  \rangle_i$.  Concretely this means that the coefficients $\tf_{ia}$ are independent Gaussian random
variables with mean 0 and variance $\eta(a)^{-1}$. In
the setting of arbitrary sparse systems [Roj96] presents a definition
and motivation of these variances that is geometric and general, in
the sense that it pertains to any sparse system.  Let $\mu_i$ be the
probability measure on $\cH_i$ that is the distribution of $\tf_i$,
and let $$\mu := \mu_1 \times \ldots \times \mu_n$$ be the
distribution of $\tf := (\tf_1, \ldots, \tf_n)$.

In the calculations used to prove Theorem 1 we also consider the model
in which the coefficient vectors $f_1, \ldots, f_n$ are statistically
independent, with each $f_i$ uniformly distributed in the unit sphere
(relative to $\|\cdot\|_i$) of $\cH_i$.  The distribution of roots
depends only on the distribution of the normalized coefficient vectors
$f_i/\|f_i\|_i$, so standard facts concerning the multivariate normal
distribution imply that, from our point of view, the two models are
equivalent.

\vskip .1cm
\goodbreak

\medskip
\Section{The Central Formula}

\medskip
We count roots in the $k$-fold product of projective spaces
$$P := P_1 \times\ldots \times P_k$$ where, for $j = 1, \ldots, k$,
$P_j := {\bf P}^{n_j}(\Re)$ is $n_j$-dimensional real projective
space.  In the usual way, the equation $f_i(\zeta) = 0$ is meaningful
for $f_i \in \cH_i$ and $\zeta \in P$ even though $f_i$ is not a
function defined on $P$. Our central concern is the expected 
number of roots $$\IndexGoal := 
\Expect{}{
\#(\{\, \zeta \in P : \tf(\zeta) = 0 \,\})},$$
but in fact we completely characterize the distribution of roots.

Let $\tZ$ be a random $n \times n$ matrix  with rows indexed
by the integers $i = 1, \ldots, n$, columns indexed by
the pairs $jh$ for $j = 1, \ldots, k$ and $h = 1, \ldots, n_j$, and
entries  $\tz_i^{jh}$ that are independently distributed normal random variables with mean zero and
variance $\delta_{ij}$.
Let $\Gamma(s) := \int_0^{\infty} {\rm exp}(-t)t^{s-1}\, dt$ be
Euler's function.

\medskip
\noindent{\bf Theorem 1:}   {\sl
\itemitem{(a)}
$$\IndexGoal = 2^{-n/2} \cdot
\Gbigl( \prod_{j = 1}^k { \Gamma({1 \over 2}) \over \Gamma({n_j
+ 1 \over 2}) } \Gbigr)
\cdot
\Expect{}{ | \det \tZ | }. 
 \myeqno \eqlabel{\MainResult}
$$
\itemitem{(b)} The induced distribution of roots is uniform: for any open $W
\subset P$,
$$\Expect{}{
\#(\{\, \zeta \in W : \tf(\zeta) = 0 \,\})} = {\vol(W) \over
\vol(P)}\IndexGoal.$$  
}

\medskip \noindent
This will be proved in Sections 8--11.  The next three sections
describe the consequences of this result.

\vskip .1cm
\goodbreak

\medskip
\Section{Reduction to Closed Form}

\medskip
In certain circumstances the RHS of (\MainResult) can be reexpressed
in closed form or in terms of the expressions derived from application
of this formula to systems that are, in certain senses, smaller.
Insofar as $\Gamma({1 \over 2}) = \sqrt{\pi}$, $\Gamma(1) = 1$, and
$\Gamma(s + 1) = s\Gamma(s)$ for all $s > 0$, the evaluations of
$\Gamma$ in (\MainResult) will be regarded as being in closed form
already, so the problem is to reduce the term $\Expect{}{ | \det \tZ |
}$.

We begin by considering systems in which there is a subset of the
variables that are determined by equations involving only those
variables.  Specifically, suppose there is some integer $k'$ between
$1$ and $k$ such that $\delta_{ij} = 0$ for all $i,j$ such that $q(i)
\le k'$ and $k' < j$,
where $q(\cdot)$ is the function defined at the end of Section 2.  Set
$n' := n_1 + \ldots + n_{k'}$.  Then $$\delta = \left[
\matrix{ \delta^{11} & 0 \cr
         \delta^{21} & \delta^{22} \cr }
\right].$$
where $\delta^{11}$, $\delta^{21}$, and $\delta^{22}$ have dimensions
$n' \times k'$, $(n - n') \times k'$, and $(n - n') \times (k - k')$
respectively.  Then (with probability one) $\tZ$ has an $n' \times (n
- n')$ block of zeros in its upper right corner, so its determinant is
the product of the determinants of the $n' \times n'$ submatrix in the
upper left and the $(n - n') \times (n - n')$ submatrix in the lower
right.  In particular, $\Expect{}{|\det \tZ|}$ does not depend on
$\delta^{21}$.  Consequently (\MainResult) implies that $\IndexGoal$
is also independent of $\delta^{21}$.  When we set $\delta^{21} = 0$
we have a cartesian product of two independent systems, and our
assumed distribution of coefficients for the combined system is the
product measure of the assumed distributions for the subsystems.  For
any particular coefficient vector for the combined system, the number
of roots is the product of the numbers of roots of the subsystems, so
the following is a consequence of the fact that the mean of a product
of independent random variables is the product of their means.
Computationally, it follows immediately from the fact that the
determinant of $\tZ$ is the product of the determinants of the
submatrices.

\medskip
\noindent {\bf Corollary 1:} {\sl
Suppose there is some $1 \le k' < k$ such that $\delta_{ij} = 0$
whenever $q(i) \le k' < j$, and let $\delta^{11}$ and $\delta^{22}$ be
as above.  Then $$E(\bn, \delta) = E((n_1,\ldots,n_{k'}),
\delta^{11}) \cdot E((n_{k' + 1},\ldots,n_{k}), \delta^{22}).$$
}

\medskip
A second general principle results from the effect on the determinant
of multiplying a row or a column by a scalar.

\medskip
\noindent{\bf Corollary 2:} {\sl
If there are nonnegative integers $d_1, \ldots, d_n$ and $e_1,
\ldots, e_k$ such that $\delta'_{ij} = d_i \cdot
e_j \cdot \delta_{ij}$, then $$E(\bn, \delta') =
\sqrt{\prod_{i = 1}^n d_i} \cdot \sqrt{\prod_{j=1}^k e_j^{n_j}}
\cdot \IndexGoal.$$  
}

\medskip
Consider now the particular case of $k = 1$ and $\delta_{11} = \ldots
= \delta_{n1} = 1$.  This corresponds to a system of $n$ linear
functionals in $n + 1$ variables, and there is exactly one projective
root for almost all coefficient vectors.  In view of (\MainResult) we
must have:

\medskip
\MMProp{ \label{\MeanAbsoluteDeterminant}
The mean absolute value of the determinant of a random $n \times n$
matrix whose entries are independently distributed normal random
variables with mean zero and unit variance is $$2^{n/2} \cdot {\Gamma({n+1 \over
2}) \over \Gamma({1 \over 2})}.$$
}

Combining the last two results with Theorem 1 yields

\medskip \noindent
{\bf Theorem 2:} {\sl
If there are nonnegative integers $d_1, \ldots, d_n$ and $e_1,
\ldots, e_k$ such that $\delta_{ij} = d_i \cdot
e_j$, then $$E(\bn, \delta) = { \Gamma({n+1 \over 2}) \over \Gamma({1
\over 2})} \cdot 
\Gbigl( \prod_{j = 1}^k { \Gamma({1 \over 2}) \over \Gamma({n_j
+ 1 \over 2}) } \Gbigr)
\cdot
\sqrt{\prod_{i = 1}^n d_i} \cdot \sqrt{\prod_{j=1}^k e_j^{n_j}}.$$  
}

\medskip \noindent
The Shub-Smale formula is the special case $k = 1$, and Rojas' formula
for unmixed systems is obtained when $d_1 = \ldots = d_n = 1$.

There is a class of systems for which $\IndexGoal$ can be computed
exactly by combining Corollaries 1 and 2 with Proposition
\MeanAbsoluteDeterminant. I know of no case outside this class 
in which the expectation $\Expect{}{ | \det \tZ | }$ evaluates to a closed
form expression.  For the systems arising from normal form games we
are able to evaluate in closed form only when $k = 2$, which
corresponds to a game with two players.  Applying ideas similar to
those underlying Corollary 1 yields:

\medskip
\noindent{\bf Corollary 3:}  {\sl
In the case of the game equilibrium
system given by (\NashEquilibrium), if $k = 2$ then
$$E(\bn, \delta) = \cases{
1  &  if $n_1 = n_2$, \cr
   \noalign{\smallskip}
0  &  otherwise.      \cr
}$$
}

\vskip .1cm

\goodbreak
\medskip
\Section{The BKK Bound for Multihomogeneous Systems}

\medskip
This section explains the consequences of Bernshtein's [Ber75] theorem
for multihomogeneous systems.  Let $f(\bz) = (f_1(\bz), \ldots,
f_n(\bz))$ be a general sparse system of $n$ equations in the $n$
variables $z_1, \ldots, z_n$, where $f_i$ has support $\cA_i \subset
{\bf N}^{n}$.  The {\it Newton polytope\/} of $f_i$ is the convex
polytope $Q_i = {\rm con}(\cA_i)$.  The {\it mixed volume\/} of $Q_1,
\ldots, Q_n$, which was first defined and studied by Minkowski, and
which we denote by ${\cal MV}( Q_1, \ldots, Q_n)$, may be defined to
be the coefficient of the monomial $\lambda_1
\cdot \ldots \cdot \lambda_n$ in the polynomial\note{See [Ewa96] for a
proof that $\vol(Q_{\lambda})$ is, in fact, a polynomial function of
$\lambda$.}  $\vol(Q_{\lambda})$ where $$Q_{\lambda} = \lambda_1Q_1 +
\ldots + \lambda_nQ_n.$$

\medskip\noindent
{\bf Theorem:} ([Ber75]) {\sl
Let $\Cstar := {\bf C} \setminus
\{0\}$.  Let $\cH^{\bf C} = \cH^{\bf C}_1 \times \ldots \times \cH^{\bf
C}_n$ where $\cH^{\bf C}_i = {\bf C}^{\cA_i}$ is the space of complex
polynomials with support $\cA_i$.  For systems $f$ in the complement,
in $\cH^{\bf C}$, of an algebraic set of positive (complex)
codimension, there are ${\cal MV}( Q_1, \ldots, Q_n)$ roots in
$(\Cstar)^n$.
}

\medskip \noindent
The maximal number ${\cal MV}( Q_1, \ldots, Q_n)$ of roots is often
referred to as the ``BKK bound'' of the system in recognition of
closely related work [Kus76, Kov78].

We apply this result to the ``demultihomogenized'' system obtained,
from the given multihomogeneous system, by setting $y_{10} = \ldots =
y_{k0} = 1$.  In comparing the roots of the latter system, in
$(\Cstar)^n$, with the roots, in $P$, of the given multihomogeneous
system, there is the possibility of roots in one of the coordinate
subspaces (in the projective sense) along which one of the variables
vanishes, but invariance under the action of $G$ quickly implies that
generic systems do not have such roots, or roots at projective
infinity.  Thus, generically, there is a one-to-one correspondence
between the roots of the given multihomogeneous system and of the
demultihomogenized system.
The Newton polytope of the $i^{\rm th}$ demultihomogenized equation is
$Q_i = \prod_{j:n_j>0} \delta_{ij}\Delj$,
where
$$\Delj := \{\, (z_{j1},
\ldots , z_{jn_j}) \in \Re^{n_j}_{\ge0} : z_{j1} + \ldots + z_{jn_j} \le 1 \,\},$$
and the generic number of complex roots of the system is
$$\BKK (\bn, \delta) := {\cal MV}\Big(\prod_{j:n_j>0} \delta_{1j}\Delj,
\ldots, \prod_{j:n_j>0} \delta_{nj}\Delj\Big).$$

Our analysis of this quantity employs the following concept. The {\it permanent\/} (e.g., [Ego96]) of an $m \times n$ matrix $D$ with
entries $d_{ij}$ is $${\rm per} \, D := \sum_{\sigma \in S_{m,n}}
\left(\prod_{i=1}^n d_{i\sigma(i)}\right)$$ where $S_{m,n}$ is the set of one to one functions from $\{1,\ldots,m\}$ to $\{1,\ldots,n\}$.  Since there are no such functions when 
$m > n$, in which case ${\rm per} \, D$ is automatically zero.  Note that multiplying any row of $D$ by a scalar has the effect of multiplying the permanent by that scalar, and that we may expand by minors along any row: for each $i = 1, \ldots,m$
$${\rm per} \, D = \sum_{j=1}^n d_{ij} \cdot {\rm per} \, D^{ij}$$
where $D^{ij}$ is the $(m - 1) \times (n - 1)$ matrix obtained from $D$ by eliminating the $i^{\rm th}$ row and the $j^{\rm th}$ column.  When $m = n$ the permanent of $D$ agrees with the permanent of its transpose, and these comments hold with rows and columns reversed. 

Let $\Delta(\bn,\delta)$ be the $n
\times n$ matrix whose first $n_1$ columns are the first column of
$\delta$, whose next $n_2$ columns are the second column of $\delta$,
and so forth.  The computation
$$\vol\Gbigl( \sum_{i=1}^n \lambda_i \Big(\prod_{j:n_j > 0}
\delta_{ij} \Delta(n_j) \Big)\Gbigr) = \vol\Gbigl(\prod_{j:n_j > 0}
\Big(\sum_{i=1}^n \lambda_i \delta_{ij}  \Big)\Delta(n_j) \Gbigr)$$
$$= {\prod_{j:n_j > 0} \Big(\sum_{i=1}^n \lambda_i \delta_{ij}  \Big)^{n_j} \over n_1! \cdot \ldots \cdot n_k!}$$
has the following immediate implication: 

\medskip
\MMProp{ \label{\RecursiveMaximum} {\rm ([McL98])} 
$$\BKK
(\bn,\delta) = {{\rm per} \, \Delta(\bn,\delta) \over
n_1!\cdot\ldots\cdot n_k!}.$$ }

\medskip
The next result enumerates consequences of the elementary properties
of the permanent, applied to this result. For $i = 1, \ldots, n$ let
$\delta^{-i}$ be the $(n - 1)\times k$ matrix obtained by eliminating
the $i^{\rm th}$ row of $\delta$.  For $j = 1, \ldots, k$ let $\be_j$
be the $j^{\rm th}$ standard unit basis vector of $\Re^k$.  In the
recursive formulas below we are adopting the convention that $$\BKK
((0,\ldots,0),\delta_0) = E((0,\ldots,0),\delta_0) = 1,$$ where
$\delta_0$ is the $0 \times k$ matrix.  This means that the ``null
system'' with no variables and no equations has one root.

\medskip
\MMProp{ \label{\RecursiveMaximumCorollaries} {\rm ([McL98])} 
\itemitem{(a)} For all $i = 1, \ldots, n$,
$$\BKK (\bn,\delta) =
\sum_{j : n_j > 0} \delta_{ij} \cdot \BKK (\bn - \be_j, \delta^{-i}).
$$
\itemitem{(b)} 
For all $j = 1, \ldots, k$ such that $n_j > 0$, $$\BKK (\bn,\delta) =
{1 \over n_j}\sum_{i = 1}^n \delta_{ij} \cdot \BKK(\bn - \be_j,
\delta^{-i}).$$
\itemitem{(c)}Suppose there is some $1 \le k' < k$ such that $\delta_{ij} = 0$
whenever $q(i) \le k' < j$, and let $\delta^{11}$, $\delta^{21}$, and
$\delta^{22}$ be as in Section 4.  Then $$\BKK (\bn, \delta) =
\BKK ((n_1,\ldots,n_{k'}),
\delta^{11}) \cdot \BKK ((n_{k' + 1},\ldots,n_{k}), \delta^{22}).$$
\itemitem{(d)} If there are nonnegative integers $d_1, \ldots, d_n$ and $e_1,
\ldots, e_k$ such that $\delta'_{ij} = d_i \cdot
e_j \cdot \delta_{ij}$, then $$\BKK (\bn, \delta') =
\Big( \prod_{i = 1}^n d_i \Big) \cdot \Big( \prod_{j:n_j>0}
e_j^{n_j} \Big)
\cdot \BKK (\bn,\delta).$$  
}

\medskip
The recursive formulas (a) and (b) give obvious algorithms for
computing $\BKK (\bn, \delta)$ that have computed values of $\BKK $ on
the order of $10^{21}$.  (Cf.~[MM97].)

In preparation for Theorem 3, we ask when
$\BKK(\bn,\delta)$ can be computed by repeated applications of (a) in
which the RHS has only one nonzero term.  We say that the pair
$(\bn,\delta)$ is {\it simply reducible\/} if the following inductive
definition is satisfied: there is some $i$ for which there is at most
one $j$ with $n_j > 0$, $\delta_{ij} > 0$, and $\BKK(\bn - \be_j,
\delta^{-i}) > 0$, and if $n > 1$ we require that for this $j$, $(\bn -
\be_j, \delta^{-i})$ is also simply reducible.  This will clearly be
the case when repeated applications of (c) reduces $\BKK(\bn,\delta)$
to a product of instances of the general homogeneous system.  In fact
this is the only way that $(\bn,\delta)$ can be simply reducible, as
we shall see in the next section.

We will need the following technical result.  Let $A$ be an $m \times
n$ matrix of 0's and 1's.  We say that an $m \times n$ matrix $D =
(d_ij)$ is $A$-{\it sparse\/} if $d_{ij} = 0$ whenever $a_{ij} = 0$.

\medskip
\MMLemma {\label{\RankCondition} The following conditions are equivalent:
\itemitem{(i)}  there is an integer $1 \le k < m$ such that, after relabelling of rows and columns, $A$ has a $k \times (n + 1 - k)$ block of 0's.
\itemitem{(ii)}  ${\rm per}\, A = 0$;
\itemitem{(iii)}  all $A$-sparse matrices have row rank less than $m$.
}

\medskip \noindent
{\bf Proof:} Clearly (i) implies (ii).  The meaning of (ii) is that
for each one-to-one $\sigma : \{1,\ldots,m\} \to \{1,\ldots,n\}$ there
is some $i$ such that $a_{i\sigma(i)} = 0$, which implies that all
$A$-sparse matrices have no $m \times m$ submatrices of full rank, so
(ii) implies (iii).  Assuming that (iii) holds, we may assume without
loss of the generality that the first $k$ rows of $A$ are minimally
linearly independent: for a generic $A$-sparse matrix $D$ their span
agrees with the span of any $(k-1)$-element subset.  Reordering
columns, we may assume that, for generic $D$, the projection of the
span of the first $k$ rows onto the space of the first $k-1$ columns
has full rank.  Now the upper right hand $k \times (n - (k-1))$ block
of $A$ must vanish, since otherwise it is straightforward to construct
an $A$-sparse matrix $D$ whose first $k$ rows are linearly
independent.
\endproof

\vskip .1cm
\goodbreak

\medskip
\Section{The Mean Exceeds the Square Root of the Maximum}

\medskip
Let $\Delta_{1 \over 2}\!(\bn, \delta)$ be the $n \times n$ matrix whose
$(i,jh)$-entry is $\sqrt{\delta_{ij}}$.  This section establishes the
following generalization of the Shub--Smale formula.

\medskip \noindent
{\bf Theorem 3:} {\sl
$${\per \, \Delta_{1 \over 2}\!(\bn,\delta) \over n_1! \cdot \ldots \cdot
n_k!}
\ge E(\bn,\delta) 
\ge \sqrt{\BKK (\bn,\delta)}
= \sqrt{{\per \, \Delta(\bn,\delta) \over n_1! \cdot \ldots \cdot
n_k!}}.$$ These inequalities hold with equality when $(\bn,\delta)$ is
simply reducible and not otherwise.  
}

\medskip
Theorem 3 will follow by induction from the following stronger result.

\medskip \noindent
{\bf Theorem 4:} {\sl
For all $i = 1, \ldots, n$,
$$\sum_{j : n_j > 0} \sqrt{\delta_{ij}} \cdot E(\bn - \be_j,
\delta^{-i})
   \ge E(\bn,\delta) \ge
\sqrt{\sum_{j : n_j > 0} \delta_{ij} \cdot 
E(\bn - \be_j, \delta^{-i})^2}, 
$$
These inequalities hold with equality if and only if there is at most one $j$ with $\delta_{ij} > 0$ and $E(\bn - \be_j, \delta^{-i}) > 0$.
}

\medskip \noindent
{\bf Proof of Theorem 3:}
The asserted inequalities follow from an induction on $n$ that
begins with the convention that $E(\bn,\delta) =
\BKK (\bn,\delta) = 1$ when $n_1 = \ldots = n_k = 0$.  The
induction step is a matter of comparing (a) of Proposition \RecursiveMaximumCorollaries\ and the
analogous formula for $\per \, \Delta_{1 \over 2}(\bn,\delta)$ with
the inequalities in Theorem 4.  Moreover, Theorem 4 implies that either of the inequalities in Theorem 3 holds with equality if and only if there is at most one $j$ with $E(\bn - \be_j, \delta^{-i}) > 0$, and $E(\bn - \be_j, \delta^{-i})$ also satisfies the inequality with equality.  In particular, it follows from induction that $ E(\bn,\delta) > 0$ if and only if
$\BKK (\bn,\delta) > 0$, so either of the inequalities in Theorem 3 holds with equality if and only if $(\bn,\delta)$ is simply reducible. \endproof

\medskip \noindent
{\bf Remark:}  We can now give a direct characterization of simple reducibility.
Applying Theorem 3 to the situation laid out in Corollary 1 and (c) of Proposition
\RecursiveMaximumCorollaries\ shows that $(\bn,\delta)$ is simply reducible if
and only if both $((n_1,\ldots,n_{k'}), \delta^{11})$ and $((n_{k' +
1},\ldots,n_{k}), \delta^{22})$ are simply reducible.  Thus it
suffices to characterize simple reducibility when the hypotheses of
(c) of Proposition \RecursiveMaximumCorollaries\ are not satisfied:
there is no $1 \le k' < k$ such that (after any reordering of rows and
columns) $\delta_{ij} = 0$ whenever $q(i) \le k' < j$.  The inequality
of Theorem 3 cannot hold with equality unless all instances of the
inequality in Theorem 4 hold with equality, so we see that if
$(\bn,\delta)$ is simply reducible, then for {\it any\/} $i = 1,
\ldots, n$ there is at most one $j$ with $n_j > 0$, $\delta_{ij} > 0$,
and $\BKK(\bn - \be_j, \delta^{-i}) > 0$, with $(\bn - \be_j,
\delta^{-i})$ simply reducible if $n > 1$.  If there is some $i$ for
which there exist distinct $j$, $j'$ with $n_j > 0$, $n_{j'} > 0$,
$\delta_{ij} > 0$, and $\delta_{ij'} > 0$, then either $\BKK(\bn -
\be_j,\delta^{-i}) = 0$ or $\BKK(\bn -
\be_{j'},\delta^{-i}) = 0$, in which case Proposition
\RecursiveMaximum\ and Lemma \RankCondition\ imply that $\Delta(\bn -
\be_j,\delta^{-i})$ or $\Delta(\bn - \be_{j'},\delta^{-i})$ has a
block of zeros, as per (iii) of Lemma \RankCondition, and this implies
that the hypotheses of (c) of Proposition
\RecursiveMaximumCorollaries\ are satisfied by $(\bn,\delta)$,
contrary to assumption.  For each $i$ there is consequently at most
one $j$ with $n_j > 0$ and $\delta_{ij} > 0$.  If there is more than
one $j$ with $n_j > 0$ it is again easy to show that the hypotheses of
(c) of Proposition \RecursiveMaximumCorollaries\ are satisfied by
$(\bn,\delta)$, so $n_j = n$ for some $j$.  That is, we have the
general homogeneous case.

\medskip
It remains to prove Theorem 4.  For the random matrix $\tZ$ of Theorem 1, let
$\tZ^i_{jh}$ be the determinant of the $(n-1) \times (n-1)$ minor
obtained by eliminating row $i$ and column $jh$.  Observe that, by
Theorem 1, $$\Expect{}{|\tZ^i_{jh}|} = 2^{n-1 \over 2}{ \Gamma({n_j
\over 2}) \over \Gamma({n_j + 1 \over 2}) } \Gbigl(
\prod_{p = 1}^k { \Gamma({n_p + 1 \over 2}) \over \Gamma({1 \over 2})
} \Gbigr)E(\bn - \be_j,\delta^{-i}),$$ so, applying Theorem 1 again to
express $E(\bn,\delta)$ in terms of $\Expect{}{|\det \tZ|}$, we
quickly find that the assertion of Theorem 4 is equivalent to:

\medskip
\MMProp{ \label{\InTermsOfZ}
For all $i = 1, \ldots, n$,
$$\sum_{j : n_j > 0} \sqrt{\delta_{ij}} { \Gamma({n_j +1
\over 2})
\over \Gamma({n_j
\over 2})} \Expect{}{|\tZ^i_{j1}|}
\ge { \Expect{}{|\det \tZ|} \over \sqrt{2} } \ge 
\sqrt{\sum_{j : n_j > 0} \delta_{ij} \Gbig( { \Gamma({n_j +1
\over 2})
\over \Gamma({n_j
\over 2})} \Gbig)^2\Expect{}{|\tZ^i_{j1}|}^2}.$$
These inequalities hold with equality if and only if $\delta_{ij}\Expect{}{|\tZ^i_{j1}|} > 0$
for at most one $j$.
}

\medskip
The proof of this will be our goal for the remainder of the section.  The next result describes the source of the inaccuracy of the approximation. 

\medskip
\MMLemma{ \label{\UpperLowerExpectedNormBounds} 
If $\tx$ is a $\Re^m_{\ge 0}$--valued random variable for which
$\Expect{}{\tx}$ is defined, then
$$\sum_{h = 1}^m \Expect{}{|\tx_h|} \ge \Expect{}{\|\tx\|} \ge
\big\|\Expect{}{\tx}\big\|.$$
The first inequality holds with equality if and only if the support of
the distribution of $\tx$ is contained in the union of the coordinate
axes. The second inequality holds with equality if and only if the
support of the distribution of $\tx$ is contained in a single ray
emanating from the origin.  }
 
\medskip \noindent
{\bf Proof:}  Since $\sum_h \Expect{}{|\tx_h|} = \Expect{}{ \sum_h  |\tx_h|}$, the first inequality follows from $\sum_{h = 1}^m  |\tx_h| \ge \|\tx\|$, and it holds
with equality if and only if, with probability one, $\sum_{h = 1}^m  |\tx_h| = \|\tx\|$.
The second inequality follows from Jensen's inequality, and it holds with equality if and only if $\|(1 - \alpha)x_0 + \alpha x_1\| = (1 - \alpha)\|x_0\| + \alpha\|x_1\|$ for any 
$x_0, x_1$ in the support of the distribution of $\tx$ and any $0 \le \alpha \le 1$.\endproof

\medskip
We will need the following technical fact. 

\medskip
\MMLemma{ \label{\MeanNormLemma}
Let $\tepsilon = (\tepsilon_1, \ldots, \tepsilon_m)$ where $\tepsilon_1,
\ldots, \tepsilon_m$ are independent identically distributed normal
random variables with mean zero and unit variance.  Then
$$\Expect{}{\|\tepsilon\|} = {\sqrt{2} \cdot
\Gamma({m+1 \over 2})\over \Gamma({m \over 2}) }. 
\myeqno \eqlabel{\MeanNorm}
$$
} 
 
\medskip \noindent
{\bf Proof:}
We compute that
$$\eqalign{
\Expect{}{\|\tepsilon\|} 
  &  = \int_{-\infty}^{\infty} \ldots
\int_{-\infty}^{\infty} \|x\| \cdot ({1 \over \sqrt{2\pi}} e^{-x_1^2/2}
\, dx_1) \cdot \ldots \cdot ({1 \over \sqrt{2\pi}} e^{-x_m^2/2}
\, dx_m) \cr
  &  = (2\pi)^{-m/2} \int_{\Re^m} \|x\| \cdot e^{-\|x\|^2/2} \, dx \cr
  &  = (2\pi)^{-m/2} \int_0^{\infty} (re^{-r^2/2}) \cdot \vol(S^{m-1})
\cdot r^{m-1} \, dr.  \cr
}$$
The asserted formula is now obtained from the
formula (e.g.,~[Fed69, p.~251])
$$\vol(S^{m-1}) = 2
{ \Gamma({1 \over 2})^m
\over
\Gamma({m \over 2}) } \quad (m \ge 1)
\myeqno \eqlabel{\SphereVolume}$$ the change of
variables $t := r^2/2$, the fact that $\Gamma({1 \over 2} ) =
\sqrt{\pi}$, and the definition of $\Gamma(\cdot)$. \endproof
 
\medskip
The next result expresses the central idea of the method, which exploits a special property of random normal variables, in its simplest form.
Random matrices
have been studied extensively [Gir90, Meh91, Mui82] but there seems to
be little prior work on mean absolute values of random determinants.

\medskip
\MMProp{ \label{\ExpansionByMinors}
Let $\tE$ be an $n \times n$ matrix whose entries $\tepsilon_{ab}$ are
independently distributed normal random variables with mean zero and
variance $\sigma_{ab}^2$.  For $1 \le a,b \le n$ let $\tE^{ab}$ be the
determinant of the $(n-1)
\times (n-1)$ minor of $\tE$ obtained by eliminating row $a$ and column
$b$.  Then for any $a = 1, \ldots, n$: 
$$
\sqrt{2/\pi} \cdot \sum_{b=1}^n \sigma_{ab} \cdot \Expect{}{|\tE^{ab}|}
\, \ge \, \Expect{}{|\det \tE|} \, \ge \,
\sqrt{2/\pi} \cdot \left( \sum_{b=1}^n \sigma_{ab}^2 \cdot \Expect{}{|\tE^{ab}|}^2
\right)^{1/2}
.$$
}

\medskip \noindent
{\bf Proof:} The expansion of the determinant by minors along 
row $a$ is $$\det \tE = \sum_{b=1}^n (-1)^{a+b} \tepsilon_{ab}\tE^{ab}.$$
For any numbers $E^{a1},\ldots,E^{an}$, elementary properties of Gaussian random variables imply that $\sum_{b=1}^n (-1)^{a+b} \tepsilon_{ab}E^{ab}$ is a normally distributed random variable with mean 0
and variance $\sum_{b=1}^n \sigma_{ab}^2 (\tE^{ab})^2.$  Since $(\tepsilon_{a1},\ldots,\tepsilon_{an})$ and $(\tE^{a1},\ldots,\tE^{an})$ are statistically independent, Fubini's theorem and 
(\MeanNorm) in the case $m = 1$ yield $$\Expect{}{|\det \tE|} = { \sqrt{2} \over
\Gamma({1 \over 2}) } {\bf E}\Bigg(\sqrt{ \sum_{b=1}^n
\sigma_{ab}^2 (\tE^{ab})^2 }  \,\,\Bigg).$$
Recalling that $\Gamma({1 \over 2}) = \sqrt{\pi}$, the claim follows
from Lemma \UpperLowerExpectedNormBounds. \endproof

Let $\Sigma_1$ be the $n \times n$ matrix with entries $\sigma_{ab}$,
and let $\Sigma_2$ be the $n \times n$ matrix with entries
$\sigma^2_{ab}$.  By an induction on $n$ we now have:

\medskip \noindent
{\bf Corollary:} {\sl
$$ 
(2/\pi)^{n/2} \cdot \per \, \Sigma_1
\ge \Expect{}{ | \det \tE | } \ge 
(2/\pi)^{n/2} \cdot \sqrt{\per \, \Sigma_2 }
.$$
}

The upper and lower bounds in Lemma \UpperLowerExpectedNormBounds\ correspond to the extreme cases in which the distribution of $\tx$ is concentrated on the coordinate axes or on the ray through $\Expect{}{\tx}$.  When the distribution of $\tx$ is known to be invariant under the action of a group, it can be possible to show that it is far from these extremes.  In the specific case we have in mind the group $$H = SO(\Re^{n_1}) \times
\ldots \times SO(\Re^{n_k})$$ acts on the space of $n
\times n$ matrices $Z$ by simultaneously acting on each row of $Z$, where each row is viewed as an element of $\Re^{n_1} \times \ldots \times
\Re^{n_k}$.  Then, because the determinant of a linear
transformation between inner product spaces is invariant under
composition with orientation preserving orthogonal transformations of
the domain or range, we have $\det(\eta Z) = \det Z$ for all $\eta \in
H$ and all $n \times n$ matrices $Z$.

Let $Z^i_{jh}$ denote the determinant of the $(n-1) \times (n-1)$
minor obtained from $Z$ by eliminating row $i$ and column $jh$.
Define the function $c_i$ from the space of $n \times n$ matrices to
$\Re^n$ by letting $c_i(Z)$ be the vector with components
$c_i^{jh}(Z) =  (-1)^{i + n_1 + \ldots + n_{j-1} + h}Z^i_{jh}$. (Of course $c_i(Z)$ is independent of the $i^{\rm th}$ row
of $Z$, and is called the {\it cross product\/} (cf.~[Spi65],
pp. 84-5) of the remaining $n-1$ rows.) 

\medskip
\MMLemma{$c_i$ is equivariant: $c_i(\eta Z) =
\eta c_i(Z)$ for all $n \times n$ matrices $Z$ and all $\eta \in H$.}

\medskip \noindent
{\bf Proof:} Let $Z_i$ denote the $i^{\rm
th}$ row of $Z$.  Then for any $\eta \in H$ we have $$Z_i \cdot c_i(Z)
= \det Z = \det(\eta Z) = \eta Z_i \cdot c_i(\eta Z).$$ Since $c_i(Z)$
and $c_i(\eta Z)$ are independent of $Z_i$, and this holds for all
$Z_i$, it must be the case that $c_i(\eta Z) =
\eta c_i(Z)$ for all $n \times n$ matrices $Z$ and all $\eta \in H$.
\endproof

\medskip \noindent
{\bf Proof of \InTermsOfZ:}  As in the last proof, 
we write $\det(\tZ) = \tZ_i \cdot c_i(\tZ)$.  As in the 
proof of Proposition \ExpansionByMinors, elementary properties of 
normal random variables and Fubini's theorem imply that 
$$\eqalign{
\Expect{}{|\det \tZ|} 
  &  = {\sqrt{2} \cdot \Gamma(1) \over \Gamma({1 \over 2})} 
{\bf E}\Bigg(\sqrt{\sum_{j=1}^k\sum_{h=1}^{n_j} \delta_{ij} \cdot
(\tZ_{jh}^i)^2}\,\,\Bigg) \cr 
}$$
$$= {\sqrt{2} \cdot \Gamma(1) \over \Gamma({1 \over 2})} 
{\bf E}\Bigg(\sqrt{\sum_{j=1}^k \delta_{ij} \cdot 
\|\txi_{ij}\|^2}\,\,\Bigg).$$ 
Combining this with Lemma \UpperLowerExpectedNormBounds\ yields
$$
{\sqrt{2} \cdot \Gamma(1)  \over  \Gamma({1 \over 2}) }\sum_{j=1}^k \sqrt{\delta_{ij}} 
\Expect{}{\|\txi_{ij}\|}
\, \ge \,
\Expect{}{|\det \tZ|}
\, \ge \,
{\sqrt{2} \cdot \Gamma(1)  \over  \Gamma({1 \over 2}) }\sqrt{\sum_{j=1}^k \delta_{ij} 
\Expect{}{\|\txi_{ij}\|}^2}. \myeqno \eqlabel{\BigInequality}$$

For $j = 1, \ldots, k$ let $\Pi_j : \Re^n \to \Re^{n_j}$ be the
projection $$\Pi_j(z^{11}, \ldots, z^{1n_1}, \,\, \ldots \,\, ,z^{k1},
\ldots, z^{kn_k}) = (z^{j1}, \ldots, z^{jn_j}).$$
Clearly $\Pi_j$ is equivariant: $\Pi_j(\eta z) =
\eta_j(\Pi_j(z))$ for all $z \in \Re^n$ and $\eta = (\eta_1, \ldots,
\eta_k) \in H$. Therefore $\Pi_j \circ c_i$ is 
equivariant.  By virtue of
elementary properties of the multivariate normal, the distribution of
the random matrix $\tZ$ on the space of $n \times n$ matrices is
invariant under the action of $H$, so the distribution of $\txi_{ij} =
\Pi_j(c_i(\tZ))$ is invariant under the action of $SO(\Re^{n_j})$.

If $\tx$ is any $\Re^{n_j}$--valued random variable whose distribution
is invariant under the action of $O(\Re^{n_j})$, the ratio
$\Expect{}{|\tx_h|}/\Expect{}{\|\tx\|}$ must agree with the mean
absolute value of the first component of a random vector that is
uniformly distributed on the unit sphere in $\Re^{n_j}$.  In
particular, by Lemma \MeanNormLemma\ we have
$${\Expect{}{\|\txi_{ij}\|} \over \Expect{}{\|\tZ_i^{j1}\|}} = {\Expect{}{\|\tepsilon\|} \over \Expect{}{|\tepsilon_1|}} = {\Gamma({n_j+1 \over 2}) /  \Gamma({n_j \over 2}) \over
\Gamma(1) /  \Gamma({1 \over 2})}$$ when $\tepsilon =
(\tepsilon_1, \ldots, \tepsilon_{n_j})$ and $\tepsilon_1, \ldots,
\tepsilon_{n_j}$ are i.i.d.~normal random
variables with mean zero.  The asserted inequality follows from substituing this into (\BigInequality).

With respect to conditions under which the inequalities hold strictly, if the vector $(\|\txi_{i1}\|,\ldots, \|\txi_{ik}\|)$ has two components 
that are nonzero with positive probability, then both inequalities in (\BigInequality)  hold strictly,
by Lemma  \UpperLowerExpectedNormBounds, and otherwise they do not. \endproof

\vskip .1cm
\goodbreak

\medskip
\Section{A Reformulation}

\medskip
This and the following three sections constitute the proof of Theorem
1.  We now reformulate the problem in a framework that is
more amenable to calculation, and which allows the application of the
methods of [SS93] and [BCS98].  Let $M_i \subset \cH_i$ be the unit sphere defined by  $\langle \cdot, \cdot \rangle_i$, and let $$M := M_1 \times \ldots
\times M_n.$$ As a submanifold of $\cH$, $M$ inherits a measure
corresponding to the intuitive notion of volume which we denote by  $\vol(\cdot)$ or (when no confusion is possible) $M$.  The {\it uniform distribution\/}
on $M$ is $\bU_M(\cdot) := \vol(\cdot)/\vol(M)$.  The analogous notation will occur in connection with other manifolds as well. The roots of $f \in \cH$ depend only on $(f_1 / \|f_1\|,
\ldots, f_n / \|f_n\|),$ and the random system $(\tf_1 / \|\tf_1\|,
\ldots,\tf_n / \|\tf_n\|)$ is uniformly distributed in $M$, by virtue of
standard facts concerning the multivariate normal distribution.

We regard $P_j$ as the space of unordered pairs $[\zeta_j] =
\{\zeta_j,-\zeta_j\}$ of antipodal points in $N_j$, where $N_j$ is the
unit sphere in $\Re^{n_j + 1}$.  Let $$N := N_1 \times \ldots \times
N_k.$$ For
each root $[\zeta] \in P$ of $f \in \cH$ there are $2^k$ corresponding roots
in $N$.

For each $i$ let $\theta_i : N \to \cH_i$ be the function with components
$\theta_{ia}(\zeta) := \eta(a)^{-1}\zeta^a$.  Let $F : M \times N \to \Re^n$ be the evaluation map with components 
$$F_i(f,\zeta) := f_i(\zeta) = \langle f_i, \theta_i(\zeta) \rangle_i.$$
The {\it incidence variety\/} is $V = F^{-1}(0)$.
Let $\pi_1$ and $\pi_2$ be the
projections from $V$ to $M$ and $N$ respectively.  We now have $$\IndexGoal = 2^{-k} \int_M \#(\pi_1^{-1}(f)) \,
d\bU_M.
\myeqno \eqlabel{\TwoTok}
$$

In preparation for the result of the next section we discuss some
technical matters.

\medskip
\MMLemma{ 
Each $\theta_i$ is equivariant with respect to the actions of $G$ on
$N$ and $\cH_i$: $\theta_i(O\zeta) = O\theta_i(\zeta)$ for all $\zeta
\in N$ and $O \in G$. The image of $\theta_i$ is contained in the unit
sphere of $\cH_i$.}

\medskip \noindent
{\bf Proof:} 
We have
$$\langle Of_i,
O\theta_i(\zeta) \rangle_i = \langle f_i, \theta_i(\zeta) \rangle_i =
f_i(\zeta) = f_i(O^{-1}(O\zeta)) = \langle Of_i, \theta_i(O\zeta)
\rangle_i.$$ Here the first equality is the invariance established in
Lemma \UniqueInvariance, and the other three equalities are
essentially matters of definition.  For given $\zeta$ this holds for all $f_i$,
so $\theta_i(O\zeta) = O\theta_i(\zeta)$.  Consequently $\|\theta_i(O\zeta)\| = \|\theta_i(\zeta)\|$ for all $\zeta$ and $O$.  Clearly $\theta_i(\zeta)$ is a standard basis vector of $\cH_i$ if $\zeta_1,\ldots,\zeta_k$ are all standard basis vectors in $\Re^{n_1 + 1}, \ldots, \Re^{n_k + 1}$ respectively, so the second claim follows from the fact that 
the action of $G$ on $N$ is transitive. \endproof

\medskip
The equation  $f_i(\zeta) = 0$ means precisely that $f_i$ and $\theta_i(\zeta)$ are orthogonal, so  for $(f,\zeta) \in V$
we may construe $\theta_i(\zeta)$ as a tangent vector in $T_{f_i}M_i$, and clearly $${\partial F_{i} \over \partial f}(f,\zeta)(0,\ldots,\theta_{i'}(\zeta), \ldots, 0)$$ is nonzero according to whether $i' = i$.  Thus $(f,\zeta)$ is a regular point of $F$, and $0$ is a regular value of $F$,  so the regular value theorem (e.g., [GP65]) implies:

\medskip
\MMLemma{ $V$ is a $C^{\infty}$ submanifold of $M \times N$ with $\dim V
= \dim M$.
}

\medskip
\def\hV{{V}}
Abusing notation, we let $V_{\zeta}$ denote both of the ``fibers''
$$\pi_2^{-1}(\zeta) \subset V \subset M \times N \quad \hbox{and } \quad 
\{\, f \in M : (f,\zeta) \in \pi_2^{-1}(\zeta) \,\}$$ over a
point $\zeta \in N$, with the appropriate interpretation to be
inferred from context.  For each $i$ let $\hV_{\zeta,i}$ be the set of
$f_i \in M_i$ with $f_i(\zeta) = 0$.  As the intersection of $M_i$
with a hyperplane, this set is a subsphere of $M_i$ of codimension
one.  Thus $\hV_{\zeta} = \hV_{\zeta,1} \times \ldots \times
\hV_{\zeta,n}$ has a simple topology that is independent of $\zeta$,
and, as one might expect:

\medskip
\MMLemma{ $\pi_2 : V \to N$ is a $C^{\infty}$ fibration.
}  

\medskip
As usual, to argue this point in detail would be a longwinded and
mundane affair, and we shall not do so.  It is, perhaps, worth
mentioning that the ``group'' of the fibration may be taken to be the group $G$
introduced in Section 3,
and that a suitable atlas of coordinate functions\note{This
terminology, and the definition of ``fibration'' we are appealing to,
are from [Ste51,
\S2].} is given by the following maps: given $\zeta_0 \in N$, a
neighborhood $W \subset N$ of $\zeta_0$, and a $C^{\infty}$ map $h : W
\to G$ satisfying $h(\zeta)\zeta_0 = \zeta$ for all $\zeta \in W$,
let $\phi: \hV_{\zeta_0} \times W \to \pi_2^{-1} (W)$ be given by
$\phi \big ( f, \zeta \big ) := \big (h (\zeta) f, \zeta \big )$. 

\vskip .1cm
\goodbreak

\medskip
\Section{An Integral Formula}

\medskip
Sard's theorem implies that almost all
points of $M$ are regular values of $\pi_1$, so we need only consider such points in computing the average number of roots.  Consider a regular point
$(f,\zeta) $ of $\pi_1$. Since $T_{(f,\zeta)}V$ is mapped surjectively onto $T_fM$ by $D\pi_1(f,\zeta)$, the restriction of $DF(f,\zeta)$ to $T_{\zeta}N \subset T_{(f,\zeta)}(M \times N)$ must be nonsingular, else $(f,\zeta)$ would not be a regular point of $F$.  The implicit function theorem implies that  there is a neighborhood $U \subset M$ of $f$ for which there is a smooth $G : U \to N$ with $G(f) = \zeta$ whose graph is contained in $V$.    The {\it condition matrix\/} at $(f,\zeta)$ is the matrix of $DG(f)$ which, by the implicit function theorem, is
$$C(f,\zeta) := - \Bigl( { \partial F \over \partial \zeta}(f,\zeta) \Bigr)^{-1}{\partial F \over \partial f}(f,\zeta) : T_fM \to T_{\zeta}N.$$  This linear transformation gives a description of the way polynomial systems $f$ are associated with roots near $(f,\zeta)$.  Let  $C^*(f,\zeta) : T_{\zeta}N \to T_fM$ be the adjoint of $C(f,\zeta)$.

\medskip \noindent
\MMProp {{\rm ([BCS98, p.~240])} \sl For any open $U \subset V$,
$$\int_M \#(\pi_1^{-1}(f) \cap U)\,dM = \int_N \int_{V_y \cap U} \det\bigl(C(f,\zeta) C^*(f,\zeta)\bigr)^{-1/2} dV_{\zeta}dN.$$} 

\medskip
\MMLemma{ If $(f,\zeta) \in V$ is a regular point of $\pi_1$, then
$$\det\big(C(f,\zeta)C^*(f,\zeta)\big)^{-1/2} = |\det Df(\zeta)|. $$
}

\medskip \noindent
{\bf Proof:}  For $v \in T_0\Re^n$ and $\phi \in T_fM = T_{f_1}M_1 \times \ldots \times T_{f_n}M_n$ we compute that
$$\eqalign{
\Bigl\langle {\partial F \over \partial f}
    (f,\zeta)\phi, v \Bigr\rangle
   &  = \Bigl\langle \bigl(\langle \phi_1, \theta_1(\zeta) \rangle_1,\ldots, \langle \phi_n, \theta_n(\zeta) \rangle_n\bigr), v \Bigr\rangle \cr
    &  = \sum_{i=1}^n \langle \phi_i, v_i\theta_i(\zeta) \rangle_i 
      = \Bigl\langle \phi, \bigl(v_1\theta_1(\zeta),\ldots,v_n\theta_n(\zeta)\bigr) \Bigr\rangle. \cr
}$$
This means precisely that the map $v \mapsto  (v_1\theta_1(\zeta),\ldots,v_n\theta_n(\zeta))$ is the adjoint ${\partial F \over \partial f}(f,\zeta)^*$ of ${\partial F \over \partial f}(f,\zeta)$, and in particular ${\partial F \over \partial f}(f,\zeta) {\partial F \over \partial f}(f,\zeta)^*$ is the identity on $T_0\Re^n$.  Since the matrix of the adjoint of a linear transformation is the transpose of the transformation's matrix, substituting the definition of the condition matrix leads to
$$\eqalign{
\det\big(C(f,\zeta)C^*(f,
  &  \zeta)\big)^{-1/2} = \Biggl(\det\Bigl( { \partial F \over \partial \zeta}(f,\zeta)^{-1} \bigl( { \partial F \over \partial \zeta}(f,\zeta) ^{-1} \bigr)^* \Bigr) \Biggr)^{-1/2}\cr
  &   = \bigl|\det { \partial F \over \partial \zeta}(f,\zeta)\bigr| = |\det Df(\zeta)|. \enspace \bull \cr
}$$

\medskip
Combining the last two results, for any open $U \subset V$ we have
$$\int_{f \in \pi_1(U)} \#(\pi_1^{-1}(f)) \, dM= \int_N \int_{V_{\zeta} \cap U}  |\det Df(\zeta)| \,\,dV_{\zeta}dN. \myeqno \eqlabel{\ToDf}$$

\vskip .1cm
\goodbreak

\medskip
\Section{Invariance}

\medskip
Combining the actions of $G$ on the various $\cH_i$ (recall Section 3) we obtain an
action of $G$ on $\cH$ given by 
$$Of := (f_1 \circ O^{-1}, \ldots, f_n
\circ O^{-1}).$$  We will exploit this symmetry to further simplify 
the RHS of the formula above.  

Each $M_i$ is invariant under the action of $G$ on $\cH_i$, of course,
so $M$ is an invariant of the action of $G$ on $\cH$, and the
restriction of this action to $M$ is an action of $G$ on $M$.  Of
course $N$ is invariant under the usual action of $G$ on
$\prod_{j=1}^k \Re^{n_j+1}$.  Combining these actions, we derive an
action of $G$ on $M \times N$ given by $O(f,\zeta) := (Of,O\zeta)$.
For any $O \in G$, $f \in M$, and $\zeta \in N$ we have $Of(O\zeta) =
f \circ O^{-1}(O\zeta) = f(\zeta)$, so:

\medskip
\MMLemma{ $V$ is an invariant of the action of $G$ on $M \times N$: 
$OV = V$ for all $O \in G$.  Consequently (for either interpretation
of the symbol $V_{\zeta}$) $O(V_{\zeta}) = V_{O\zeta}$ for all $\zeta$
and $O$.  }

\medskip \noindent
\MMProp{
The quantity $\int_{\hV_{\zeta}} |\det Df(\zeta)| \, dV_{\zeta}$ is
independent of $\zeta$.
}

\medskip \noindent
{\bf Proof:} Observe that $$D(Of)(O\zeta) = D(f \circ O^{-1})(O\zeta)
= Df(\zeta) \circ O^{-1}$$ so that $ |\det D(Of)(O\zeta)| = |\det
Df(\zeta)|.$ We now have the calculation that $$\int_{\hV_{O\zeta}}
|\det Df(O\zeta)| \, d\hV_{O\zeta} =
\int_{\hV_{\zeta}} |\det D(Of)(O\zeta)| \, d\hV_{\zeta} = \int_{\hV_{\zeta}}
|\det Df(\zeta)| \, d\hV_{\zeta}.$$ Here the first equality is an application of
the change of variables formula with the change of variables function
an isometry, so that the Jacobean is identically one.  The claim now
follows from the fact that the action of $G$ on $N$ is transitive.
\endproof

Applying this to (\ToDf), for any open $W \subset N$ and any $\zeta
\in N$ we have $$\int_{M}\#(\pi_1^{-1}(f) \cap \pi_2^{-1}(W)) \, dM = \vol(W)
\cdot \int_{\hV_{\zeta}} |\det Df(\zeta)| \, d\hV_{\zeta}. 
\myeqno \eqlabel{\AfterInvariance}$$ 
Clearly (b) of Theorem 1 follows directly from this.  The remaining
task is to prove (a) of that result.

%\medskip
%\MMProp{ \label{\AfterInvariance}
% For any open $W \subset N$ and any $\zeta \in N$,
%$$\int_{M}\#(\pi_1^{-1}(f) \cap \pi_2^{-1}(W)) \, dM = \vol(W)
%\cdot \int_{\hV_{\zeta}} |\det Df(\zeta)| \, d\hV_{\zeta}.$$ 
%}

%\medskip
%Clearly (b) of Theorem 1 follows directly from this.  The remaining
%task is to prove (a) of that result.
%\endprop

\vskip .1cm
\goodbreak

\medskip
\Section{The Final Calculations}

\medskip
Fixing $\zeta \in N$, let  $\tf_{\zeta} = (\tf_{\zeta,1}, \ldots,
\tf_{\zeta,n})$ be the orthogonal projection of $\tf$ onto the
subspace of polynomial systems for which $\zeta$ is a root.
For each $i$, $\|\tf_{\zeta,i}\|$ and
$\tf_{\zeta,i}/\|\tf_{\zeta,i}\|$ are statistically independent, and
the normalized vector is uniformly distributed in $V_{\zeta,i}$, so
$$\eqalign{
\int_{\cH}  |\det D\tf_{\zeta}(\zeta)| \, d\mu 
  &  = \int_{\cH} \Gbig(\prod_{i=1}^n \|\tf_{\zeta,i}\|\Gbig) 
\cdot  \Big|\det D\Big({ \tf_{\zeta,1} \over \|\tf_{\zeta,1}\| }, 
\ldots, { \tf_{\zeta,n} \over \|\tf_{\zeta,n}\| }\Big) (\zeta)\Big| 
\, d\mu \cr
  &  = \Gbig( \prod_{i=1}^n \Expect{}{\|\tf_{\zeta,i}\|} \Gbig) 
\int_{\hV_{\zeta}} |\det Df(\zeta)| \, d\bU_{\hV_{\zeta}}. \cr
}$$ Combining this with (\TwoTok) and (\AfterInvariance), we now
obtain $$E(\bn,\delta) = 2^{-k}{ \vol(N) \cdot \vol(V_{\zeta}) \over
\vol(M) \cdot \prod_{i=1}^n \Expect{}{\|\tf_{\zeta,i}\|} } 
\int_{\cH}  |\det D\tf_{\zeta}(\zeta)| \, d\mu.$$

The formula (\SphereVolume) for sphere volume gives $$\vol(N_j) = 2{
\Gamma({1 \over 2})^{n_j + 1} \over \Gamma({n_j +1 \over 2}) }, \quad
\vol(M_i) = 2{ \Gamma({1 \over 2})^{\dim \cH_i} \over 
\Gamma({\dim \cH_i \over 2}) },
\quad \vol(V_{\zeta,i}) = 2{ \Gamma({1 \over 2})^{\dim \cH_i - 1}
\over \Gamma({\dim \cH_i - 1 \over 2}) },$$ and Lemma \MeanNormLemma\
yields $$\Expect{}{\|\tf_{\zeta,i}\|} = { \sqrt{2} \cdot
\Gamma({\dim \cH_i \over 2}) \over \Gamma({\dim \cH_i - 1 \over 2})
}.$$ Since $\vol(M) = \vol(M_1) \times \ldots \times \vol(M_n)$, and
similarly for $N$ and $V_{\zeta}$, we now have
$$\IndexGoal = 2^{-n/2} \cdot
\Gbigl( \prod_{j = 1}^k { \Gamma({1 \over 2}) \over \Gamma({n_j
+ 1 \over 2}) } \Gbigr)
\cdot
\int_{\cH}  |\det D\tf_{\zeta}(\zeta)| \, d\mu. 
\myeqno \eqlabel{\NextToLast}$$

In the further evaluation of this quantity we are free to let $\zeta$
be any convenient point in $N$. We will compute at $\zeta_0 =
(\be_{10}, \ldots, \be_{k0}) \in N$ where, for $1 \le j \le k$,
$\be_{j0}, \be_{j1}, \ldots,
\be_{jn_j}$ are the standard unit basis vectors of $\Re^{n_j + 1}$.  
For each $i$ and $j$ let $a^0_{ij} = (\delta_{ij},0,\ldots,0) \in \cA_{ij}$,
and for each $i$ let $a_i^0 = (a^0_{i1}, \ldots, a^0_{ik}) \in \cA_i$.  Since
$\zeta_0^a = 0$ for all $a \in
\cA_i$ other than $a^0_i$, and $\zeta_0^{a^0_i} = 1$,  for each $i$
$ \hV_{\zeta_0,i} := \{\, f_i \in M_i : f_{ia_i^0} = 0 \,\}.$
For $i$ and $j$ such that $\delta_{ij} > 0$ and each $h = 1, \ldots,
n_j$, let $a^{jh}_i$ be $a^0_i$ with $a^0_{ij}$ replaced by
$(\delta_{ij} - 1, 0,\ldots,0,1,0,\ldots,0)$ (the `1' is component
$h$).  Then $T_{\zeta_0}N$ is spanned by the $n$ vectors
$${\bf b}_{jh} := (0, \ldots, \be_{jh},
\ldots,0) \qquad (1 \le j \le k, 1
\le h \le n_j),$$ and elementary calculus yields 
$$D\tf_{\zeta_0,i}(\zeta_0){\bf b}_{jh} = \cases{ \tf_{ia_i^{jh}} & if $\delta_{ij}
> 0$, \cr 0 & if $\delta_{ij} = 0$. \cr}$$ In this way we obtain a
description of $D\tf_{\zeta_0}(\zeta_0)$ as an $n \times n$ matrix with rows
indexed by $f_1,\ldots, f_n$, columns indexed by the pairs $(j,h)$,
and this $(i,jh)$--entry.  Recalling from Section 3 that the variance of $
\tf_{ia_i^{jh}}$ is $\delta_{ij}$, we see that the matrix of $D\tf_{\zeta_0}(\zeta_0)$
has the same distribution as $\tZ$.  In view of (\NextToLast) this observation
completes the proof of Theorem 1.

\bigskip
\vfill \eject
\goodbreak

% \vfill \eject
% \newcount\bibnumber
% \def\bibitem{\global\advance\bibnumber by 1
% 	\hangindent=36pt \the\bibnumber. }

% \def\bibitem{\global\advance\bibnumber by 1
% 	\itemitem{\the\bibnumber.}}

\newcount\bibnumber
\def\bibitem#1#2{\smallskip \hangindent=36pt $\hbox to 1.4cm{#1 \hfil}$#2}
\bibnumber = 0

\parindent=0pt
\centerline{\bf References}

\medskip
\bibitem{[Ber75]} {D.~N.~Bernshtein, The number 
of roots of a system of equations, {\it Functional Analysis and its
Applications} {\bf 9} (1975), 183--185.}

\bibitem{[BCS93]} { L.~Blum, F.~Cucker, M.~Shub, and S.~Smale, {\it Complexity
and Real Computation\/},  Springer-Verlag, New York, (1998).}

\bibitem{[BP32]} {A.~Bloch and G.~P\'olya, On the roots of a certain algebraic equation, {\it Proc.~London Math.~Soc\/}.~ {\bf 33} (1932), 102--114.}

%\bibitem{[CS61]} {R.~D.~Carmichael and E.~R.~Smith, {\it Mathematical 
%Tables and Formulas\/}, Dover, New York, (1961).}

%\bibitem{[CE95]} {J.~F.~Canny and I.~Emiris, Efficient Incremental
%Algorithms for the Sparse Resultant and the Mixed Volume,
%{\it Journal of Symbolic Computation\/} {\bf 20} (1995), 117--149.}

%\bibitem{[Dre]} {M.~Dresher, Probability of a pure equilibrium
%point in $n$-person games, {\it Journal of Combinatorial Theory\/}
%{\bf 8} (1970), 134--145.}

\bibitem{[EK95]} {A.~Edelman and E.~Kostlan, How many zeros of a random
polynomial are real?, {\it Bulletin of the American Mathematical
Society\/} {\bf 32} (1995), 1--37.}

\bibitem{[Ego96]} {G.~P.~Egorychev, Van der Waerden conjecture and
applications, {\it Handbook of Algebra\/}, Vol.~I, Elsevier,
Amsterdam (1996), 3--26.}

\bibitem{[Ewa96]} {G.~Ewald, {\it Combinatorial Convexity and
Algebraic Geometry\/}, \break Springer, 
New York, (1996).}

\bibitem{[Fed69]} {H.~Federer, {\it Geometric Measure Theory\/}, Springer, 
New York, (1969).}

\bibitem{[FT91]} {D.~Fudenberg and J.~Tirole, {\it Game Theory\/}, MIT Press, 
Cambridge, (1991).}

\bibitem{[Gir90]} {V.L.~Girko, {\it The Theory of Random Determinants\/}, 
Kluwer, \break Boston, (1990).}

\bibitem{[GP65]} {V.~Guillemin and A.~Pollack, {\it Differential$\!$
Topology\/}, Prentice-Hall, Englewood Cliffs, (1965).}

% \bibitem{[GPS]}{F.~Gul, D.~Pearce, and E.~Stachetti, A
% bound on the proportion of pure strategy equilibria in generic
% games, {\it Mathematics of Operations Research\/} {\bf 18} (1993), 548--
% 552.}

%\bibitem J.~C.~Harsanyi, Oddness of the number
%of equilibrium points: a new proof, {\it Int.~J.~Game Theory}, {\bf
%2} (1973), 235--250.  

\bibitem{[Kac43]} {M.~Kac, On the average number of real roots of a
random algebraic equation, {\it Bulletin of the American Mathematical
Society\/} {\bf 49} (1943), 314--320 and 938.}  

%\bibitem H.~Keiding, On the maximal number of Nash equilibria 
%in a bimatrix game.  Forthcoming in {\it Games and Economic
%Behavior}.  

\bibitem{[Kho78]} {A.~G.~Khovanskii, Newton 
polyhedra and the genus of complete intersections, {\it Functional Analysis and its
Applications} {\bf 12} (1978), 51--61.}

%\bibitem E.~Kohlberg and J.-F.~Mertens, On the strategic
%stability of equilibria, {\it Econometrica\/} {\bf 54} (1986),
%1003--1038.  

\bibitem{[Kos93]} E.~Kostlan, On the distribution of roots of random
polynomials, {\it From Topology to Computation: Proceedings of the
Smalefest\/}, Hirsch, M., Marsden, J., and Shub, M. (eds) (1993).

%\bibitem H.~Kuhn, et.~al., The work of John Nash in game
%theory: Nobel seminar, December 8, 1994, {\it Journal of Economic
%Theory\/} {\bf 69} (1996), 153--185.  

\bibitem{[Kus75]} {A.~G.~Kushnirenko, The Newton 
polyhedron and the number of solution of a system of $k$ equations in
$k$ unknowns, {\it Upsekhi Mat. Nauk.} {\bf 30} (1975), 266--267.}

\bibitem{[MM97]} {R.~D.~McKelvey and A.~McLennan, The
maximal number of regular totally mixed Nash equilibria, {\it Journal
of Economic Theory} {\bf 72} (1997), 411--425.} 

%\bibitem A.~McLennan, The
%maximal generic number of pure Nash equilibria, {\it Journal of
%Economic Theory} {\bf 72} (1997a), 408--410.}  

\bibitem{[McL97]} {A.~McLennan, On the Expected Number of Nash Equilibria
of a Normal Form Game, mimeo, University of Minnesota, (1997).}

\bibitem{[McL98]} {A.~McLennan, The maximal number of real roots of a 
multihomogeneous system of polynomial equations, forthcoming in {\it
Beitr\"age zur Algebra und Geometrie\/}, (1998).}

%\bibitem A.~McLennan and I.-U.~Park, Generic 4x4 games have at 
%most 15 Nash equilibria, Discussion Paper 300, Center for Economic
%Research, University of Minnesota, (1997). 

\bibitem{[Meh91]} M.L.~Mehta, {\it Random Matrices}, Academic Press,
New York \break (1991). 

\bibitem{[Mui82]} R.J.~Muirhead, {\it Aspects $\!$of$\!$ Multivariate Statistical Theory}, Wiley, New York (1982). 

% \bibitem F.~Morgan, {\it Geometric Measure Theory: a Beginner's
% Guide}, Academic Press, New York (1988).

%\bibitem J.~Nash, ``Non-cooperative games,'' 
%Ph.D.~Thesis, Mathematics Department, Princeton University, (1950).

%\bibitem J.~Nash, Non-cooperative games, 
%{\it Annals of Mathematics\/} {\bf 54} (1951), 286--295. 

%\bibitem M.~Okuno-Fujiwara and A.~Postlewaite, Social Norms
%and Random Matching Games, {\it Games and Economic Behavior\/} {\bf 9}
%(1995), 79--109.  

%\bibitem I.~Powers, Limiting distributions of the number of
%pure strategy Nash equilibria, {\it International Journal of Game
%Theory\/} {\bf 19} (1990), 277--286.  

%\bibitem T.~Quint and M.~Shubik, A theorem on the number of
%Nash equilibria in a bimatrix game, {\it International Journal of Game
%Theory\/} {\bf 26} (1997), 353--360.  

\bibitem{[Roj96]} {J.~M.~Rojas, On the average number of real roots
of certain random sparse polynomial systems, {\it Lectures on Applied
Mathematics Series\/}, ed. by J.~Renegar, M.~Shub, and S.~Smale,
American Mathematical Society, (1996).} 

\bibitem{[SS93]} {M.~Shub and S.~Smale, Complexity of Bezout's
theorem II: volumes and probabilities, {\it Computational Algebraic
Geometry\/} (F. Eyssette and A.~Galligo, eds.), Progr.~Math.,
vol.~109 (1993), Birkhauser, Boston, 267--285.}

\bibitem{[Spi65]} {M.~Spivak, {\it Calculus on Manifolds : a Modern 
Approach to Classical Theorems of Advanced Calculus\/}, Benjamin, New
York, (1965). }

%\bibitem M.~Spivak, {\it  A Comprehensive Introduction to 
%Differential Geometry. Volume 1.\/}, 2nd ed., Publish or Perish,
%(1979).

%\bibitem W.~Stanford, The Poisson distribution and
%probability of $k$ pure equilibria in matrix games, mimeo, Department of
%Economics, University of Illinois at Chicago, (1993).  

% \bibitem W.~Stanford, The limit distribution of the
% number of pure equilibria in symmetric bimatrix games, mimeo,
% Department of Economics, University of Illinois at Chicago, (1994).

\bibitem{[Ste51]} {N.~Steenrod, {\it The Topology of Fibre Bundles\/},
Princeton University Press, Princeton (1951).}

\bigskip
\baselineskip10pt
\bigskip \noindent {\smallsmc
	Andrew McLennan \hfill \break \noindent 
	Department of Economics \hfill\break \noindent
	University of Minnesota\hfill\break \noindent
	271 19th Ave. S.\hfill\break \noindent
	Minneapolis, MN  55455, USA\hfill \break \noindent
\par}

\noindent {\tt mclennan@atlas.socsci.umn.edu} \hfill \break
\noindent {\tt http://www.econ.umn.edu/{\char'176}mclennan}

\end